\documentclass[11pt,a4paper]{article}
\usepackage{amssymb, amstext, amscd, amsmath, color}
\usepackage{graphicx, url, hhline}
\usepackage[left=2.50cm, right=2.00cm, top=2.30cm, bottom=2.00cm]{geometry}

\usepackage{float}

\newtheorem{defin}{}
\newtheorem{saetze}[defin]{}
\newtheorem{conjec}[defin]{}
\newtheorem{lemmas}[defin]{}
\newtheorem{folger}[defin]{}
\newtheorem{bemerk}[defin]{}

\newenvironment{theorem}  {\begin{saetze}\it {\bf Theorem:}}{\end{saetze}}

\newcommand{\fillbox}{\mbox{$\bullet$}}

\newcommand{\Z}{\mathbb Z}

\newcommand{\R}{\mathbb R}
\newcommand{\GAP}{\mathsf{GAP}}

\begin{document}

\bigskip
\bigskip

\centerline{\large \textbf{Short presentations for crystallographic groups}}

\bigskip
\bigskip

\centerline{\large Igor A. Baburin}
\centerline{\emph{Email: baburinssu@gmail.com}}

\bigskip

\noindent A practical approach is proposed to construct short presentations for Euclidean crystallographic groups in terms of generators and relations. For our purposes a short presentation is the one with a small number of short relators for a given generating set. The connection is emphasized between relators of a group presentation and cycles in the associated Cayley graph. It is shown by examples that a short presentation is usually the one where relators correspond to \emph{strong rings} in the Cayley graph and therefore provide a natural upper bound for their size. Presentations are computed for vertex-transitive groups which act with trivial vertex stabilizers on a number of high-symmetry 2-, 3- and 4-periodic graphs. Higher-dimensional as well as subperiodic examples are also considered. Relations are explored between geodesics in periodic graphs and corresponding cycles in their quotients.


\section{Introduction}

\medskip

A presentation of a group is its abstract definition in terms of generators and relations. Instead of relations it is possible to speak of \emph{relators} that are expressions in generators evaluating to the identity of a group (for example, from the relation $ab=cd$ a relator $abd^{-1}c^{-1}$ is obtained) -- the convention that we shall follow in this paper. For a group $G$ we write $G = \langle s_{1}, ..., s_{k} | R \rangle$ where $s_{1},..., s_{k}$ are its generators, and $R$ is a set of \emph{defining relators} on $\{ s_{1}, {s_{1}}^{-1}, ... , s_{k}, {s_{k}}^{-1} \}$. A group is said to be finitely generated and finitely presented if both sets of generators and relators are finite (moreover, relators are all of finite length). For example, the dihedral group of order 4 admits the following presentation on three generators: $D_{2} = \langle a, b, c | a^{2}, b^{2}, c^{2}, abc \rangle$ that can be determined from the action of $D_{2}$ on the vertices of a tetrahedron.

Group presentations play a central role in group-theoretic computations, historically the first being subgroup computations \cite{Dietze74}. In the classical book on discrete groups Coxeter and Moser~(1980) \cite{Coxeter80} dedicated the whole chapter to presentations of two-dimensional crystallographic groups (so-called \emph{wallpaper groups}) that were treated using an elementary geometrical approach. The idea of this paper is to propose a practical recipe that makes it possible to do -- what Coxeter and Moser did by hand -- fully automatically for any generating set and in any dimension.

Presentations of crystallographic groups find important applications outside pure group theory. Based on the theorem of Poincar\'e on fundamental polygons, Delone (1959) \cite{Delone59} proposed a method to describe (and enumerate) \emph{isohedral} plane tessellations by the so-called \emph{adjacency symbols} that are essentially group presentations coded in a special way. The ideas of Delone were developed further by Gr\"unbaum and Shephard \cite{Gruen77}. In the Delaney--Dress combinatorial theory of tilings an essential part consists in recognizing groups from finite presentations, in order to decide on the realizability of a tiling in a certain space \cite{Delgado2001}. In connection with Euclidean space forms Moln{\'a}r (1987) derived minimal presentations for fixed-point-free crystallographic groups (Bieberbach groups) in $\R^3$ by a suitable choice of fundamental domains \cite{Molnar87} -- the approach that he applied later to other space groups \cite{Molnar88}. Lord (2003) compiled a table with generators and relations (one for each group) for all space groups in $\R^3$ \cite{Lord2003}.

In crystallography it was understood long time ago (W. Fischer and collaborators) that the complete derivation of homogeneous sphere packings -- with a few minor exceptions \cite{FK95} -- can be restricted to general positions of space groups \cite{Fischer68, Fischer71}. In this case a (connected) sphere packing is obtained if sphere contacts correspond one-to-one to space-group generators. Therefore, enumeration of sphere packings is essentially reduced to the enumeration of generating sets for space groups with an additional property that the images of a point under group generators produce shortest equal distances with the original point (\emph{i.e.}, the center of a reference sphere). In this way subperiodic sphere packings \cite{KF78} and 3-periodic sphere packings with contact number 3 \cite{Fischer83, KF95} were derived. Complete results for homogenous 3-periodic sphere packings are available for all crystal systems except for the general position of $C2/c$. Therefore, sphere packings provide a rich collection of `natural' generating sets for three-dimensional crystallographic groups. No attempt has been made, however, to write out group presentations in terms of these `natural' sets of generators.

Presentations of crystallographic groups can be constructed using built-in functions from computer algebra packages like $\GAP$ \cite{GAP} or \textsc{Magma}. $\GAP$ contains a library of crystallographic groups in dimensions $d \leq 6$ from which their presentations can be routinely computed. The computed presentations are, however, rather lengthy (as it is the case for \emph{e.g.} polycyclic presentations). From a theoretical point of view it is no problem to rewrite a given presentation in terms of desired generators. To do so, the expressions for old generators in terms of new ones need to be put into old relators. Tietze transformations can be applied to simplify the `as-rewritten' presentation. Our experience suggests, however, that rewriting `default' presentations of crystallographic groups often produces rather long and unnecessarily complicated relators which cannot be efficiently dealt with by Tietze transformations.

In computational group theory constructing \emph{short} (also referred to as \emph{concise}, see \cite{Cannon73}) presentations for a given set of generators is a well-known area of activity. For finite groups effective algorithms exist for computing such presentations \cite{Cannon73, Strogova98}. The algorithms are based essentially on finding a special cycle basis (that is in a certain sense minimal) for a Cayley graph (see below).

The aim of this paper is to propose a practical approach for computing short presentations of crystallographic groups on given generating sets and to study the associated Cayley graphs that are interesting for crystallography and solid-state sciences. Herein, by the term \emph{short} we mean a small number of short relators for a given generating set. In many situations our approach yields a non-redundant set of shortest relators -- however, this cannot be guaranteed in all cases. Another motivation for this study is the intention to find a good (in the lucky case, the exact) upper bound for the sizes of strong rings in vertex-transitive periodic graphs. Recall that in a crystallographic terminology a \emph{strong ring} is a cycle that cannot be written as a sum (mod 2) of strictly smaller cycles \cite{Goez91}. Strong rings are usually found by increasing iteratively the size of cycles to be included in the decomposition. The exact upper bound is generally not known and could not be proven to be found unless additional arguments are invoked (from \emph{e.g.} the tiling construction, see \cite{Blatov07}). Such arguments are, however, not always available. From the connection between relators of a group presentation and cycles in the respective Cayley graph we aim at computing a good upper bound for the size of strong rings in \emph{any} graph from a large class of vertex-transitive graphs, that is the Cayley graphs of crystallographic groups, -- without any reference to an embedding or drawing.

\section{Theoretical and computational background}
\subsection{Presentations of crystallographic groups}

We consider crystallographic groups of a $d$-dimensional Euclidean space $\R^d$. From an algebraic point of view a crystallographic group $G$ is an extension of a
free abelian group $T \cong \Z^d$ of rank $d$ by a finite group $P$ so that $P$ acts faithfully on $T$. The group $T$ is the translation subgroup 
of $G$ and $P$ is its point group.

A presentation of a group extension can be constructed from a presentation of a normal subgroup and a presentation of a quotient as given by the following theorem \cite[pp.~38-39]{Handbook}.

\begin{theorem} Let $G$ be an extension of a normal subgroup $T$ by $G/T$. Assume we have presentations $\langle Y|S \rangle$ of $T$ and $\langle \bar{X}|\bar{R} \rangle$ of $G/T$ on generating sets $Y$ and $\bar{X}$, resp. For homomorphism $\phi: G \rightarrow G/T$, let $X = \phi^{-1}(\bar{X})$. We define $R = \{ rw_r^{-1} | r \in \phi^{-1}(\bar{R}) \}$ for $w_r \in (Y \cup Y^{-1})^{\ast}$ and $M = \{ x^{-1}yxw_{xy}^{-1} | x \in X, y \in Y \}$ for $w_{xy} \in (Y \cup Y^{-1})^{\ast}$. Then $\langle X \cup Y|R \cup S \cup M \rangle$ is a presentation of $G$.
\end{theorem}

\noindent This theorem carries over as is to subperiodic crystallographic groups (when $rank(T) < d$), in which case the point group $P$ need not act faithfully on the lattice $T$ \cite{Koe80}. In $\R^3$ these are known as \emph{rod groups} ($rank(T)=1$) and \emph{layer groups} ($rank(T)=2$). \\

Given a set of generators $S$ (in the form of \emph{augmented} $(d+1) \times (d+1)$ matrices) for a crystallographic group $G$, we have to go through the following steps to obtain the relators:

\begin{enumerate}
\setlength\itemsep{-0.02in}
\item[(a)] Compute a short presentation for the point group $P = G/T$ on the generators which are images of those from $S$ under homomorphism $G \rightarrow P$.
\item[(b)] Write out a presentation for $T$ in terms of $S$.
\item[(c)] Replace the generators in the relators of $P$ by their preimages from $S$. If any relator now corresponds to a non-trivial translation, multiply it by a negative translation expressed in the basis of $T$ (that is in turn given in terms of $S$).
\item[(d)] Compute conjugates of the generators of $T$ by elements from $S$ and form the relators from the set $M$.
\end{enumerate}

\noindent Step (a) can be accomplished using the approach of Cannon (1973) for finite groups \cite{Cannon73} as implemented \emph{e.g.} in the $\GAP$ function `IsomorphismFpGroupByGenerators'. Shortest expressions for the generators of $T$ in terms of $S$ are found at Step (b) by a breadth-first traversal of the respective Cayley graph. The relators of $T \cong \Z^d$ are usually chosen as $d(d-1)/2$ commutators of the $d$ free generators (expressed in terms of $S$) that are later used at Steps (c-d). Finally, Tietze transformations are applied to simplify the obtained relators.

\subsection{Cayley graphs}

To a group $G$ with an inverse-closed generating set $S \subset G \setminus \left \{1_G\right \}$ (such that $S = S^{-1}$) it is possible to associate a connected undirected graph without loops and multiple edges, a so-called \emph{Cayley graph}, whose vertices correspond one-to-one to elements of $G$. Two vertices $g, h \in G$ are joined by an edge (corresponding to generator $s$) if $h = gs$ holds for some $s \in S$. The action of $G$ on itself by left multiplication is also the action of $G$ on vertices of a Cayley graph as a group of automorphisms, \emph{i.e.,} adjacency-preserving permutations of vertices. This action of $G$ is both \emph{free} (that is, with trivial vertex stabilizers) and \emph{transitive}, and in group theory is termed \emph{regular}.\footnote{Here the term \emph{regular} is used in accordance with existing group-theoretic terminology and should not be confused with its meaning in geometry (\emph{e.g.} regular polyhedra). For physisists our usage of the term \emph{regular} should be familiar from representation theory.} For example, a square is a Cayley graph for two groups of order 4: the cyclic group $C_4$ and the dihedral group $D_2$ (note that isomorphic graphs can arise as Cayley graphs of non-isomorphic groups). Given a presentation $\langle S|R \rangle$ of $G$, any path along the edges of a Cayley graph is a word in generators from $S$. Any relator from $R$ (of length $>2$) corresponds to a cycle in a Cayley graph \cite{Coxeter80} and furthermore, relators can be chosen in a way so as to correspond to the cycles from a cycle basis \cite{Cannon73, Coxeter80}.

Cayley graphs of crystallographic groups in $\R^d$ -- that are our main focus in this paper -- are special cases of $d$-\emph{periodic graphs} which are simple undirected graphs with a finite number of vertex and edge orbits mod $\Z^d$ \cite{Delgado2004}. A $d$-periodic graph is \emph{minimal} \cite{Klee92} if it is connected and has the minimum possible number of vertex- and edge-orbits mod $\Z^d$.

Consider a vertex-transitive graph $\Gamma$ with a full automorphism group $Aut(\Gamma)$ and let $Stab_{Aut(\Gamma)}(x)$ be the stabilizer of an arbitrary vertex $x$. Then $\Gamma$ arises as a Cayley graph for any regular subgroup of $Aut(\Gamma)$. Whether or not a regular subgroup exists can be determined as follows. The index of a regular subgroup is equal to $|Stab_{Aut(\Gamma)}(x)|$. As in our practice we encounter only graphs with finite vertex stabilizers in $Aut(\Gamma)$, we first compute representatives of conjugacy classes for all subgroups of $Aut(\Gamma)$ with index = $|Stab_{Aut(\Gamma)}(x)|$, from which regular subgroups are filtered out based on the following condition \cite{Sene88}: let $H<Aut(\Gamma)$ and $K = Stab_{Aut(\Gamma)}(x) \cap H$. The subgroup~$H$ is regular iff the intersection group $K$ is trivial. The computation of subgroups can be carried out from a convenient finite presentation of $Aut(\Gamma)$ using algorithms for finitely-presented groups. In case $Aut(\Gamma)$ is polycyclic (that is the case for crystallographic groups in dimensions $d \leq 3$ and mostly in $d=4$), then algorithms for polycylic groups can be applied \cite{Eick24}. The intersection group is computed by a slightly adapted method for finite matrix groups. Furthermore, generating sets of regular subgroups are determined which consist of elements mapping vertex $x$ to the adjacent vertices, and precisely these generating sets are used in the following for computing presentations.

\subsection{Examples}

\subsubsection{Presentations for translation groups}

If translation group $T \cong \Z^d$ is generated by $n > d$ free generators, then relators can be written out in two (in principle, equivalent, but technically different) ways: either $n(n-1)/2$ commutators are combined with ($n-d$) relators or alternatively, a basis is first determined and only the $d(d-1)/2$ commutators for the basis vectors are considered whereas the remaining relators are derived from the expressions for the generators in terms of the basis found.

\noindent \emph{Example:} consider the generating set of three vectors for $\Z^2$ that corresponds to the familiar $3^{6}$ plane tiling by triangles: $\{ a=(1, 0), b=(0, 1)$, $c=(1,1) \}$. The two ways of deriving relators give rise to the following presentations: (a) $\langle a, b, c | [a, b], [b, c], [a, c], abc^{-1}\rangle$ and  (b) $\langle a, b, c | [a, b], abc^{-1}\rangle$ (square brackets as usual stand for the commutators of the respective elements: $[a,b]=a^{-1}b^{-1}ab$). Both presentations are simplified to the shortest one $\langle a, b, c | abc^{-1}, c^{-1}ba \rangle$ (two relators of length 3 correspond to two triangles per unit cell) but the effort will be more considerable in case (a), especially in higher dimensions.

Turning to a slightly more general situation, namely, the generating sets for $\Z^2$ consisting of vectors $\{ a=(1, 0), b=(0, 1), c=(n, n) \}$ ($n>1$ is an arbitrary integer), we arrive at the sequence of Cayley graphs for $\Z^2$ that converges (in the sense of \cite{Trofimov2020}) as soon as $n \rightarrow \infty$ to a three-dimensional cubic lattice. That is, each Cayley graph from this sequence is isomorphic to a cubic lattice within the ball of radius $n$. The presentation for $\Z^2$ in this case clearly reads as: $\Z^2 = \langle a, b, c | [a, b], (ab)^{n}c^{-1}\rangle$. Note the odd cycles of size $2n+1$ meaning that starting from the $(n+1)^{th}$ sphere onwards the members of the sequence start to be distinguished from a cubic lattice.

Non-minimal generating sets for $\Z^d$ become relevant for when shortest words correspond to translations that do not form a lattice basis. In the {\bf elv} net (the RCSR database of M. O'Keeffe) that is a Cayley graph for the $I23$ space group with the generating set \\
\centerline{$\{ a=x, -y, -z; b=-x, y, 1-z; c=1/2+z, 1/2-x, 3/2-y \}$} \\
shortest words corresponding to translations are $c^3 = (3/2, -3/2, 3/2), (ab)^2 = (0, 0, 2), acbc^{-1} = (0, 2, -1), bc^{-1}ac = (1, 0, 2)$ which altogether (but not any three of them) generate the full $cI$-lattice whereas the primitive translations $(-1/2, 1/2, 1/2), (1/2, -1/2, 1/2), (1/2, 1/2, -1/2)$ correspond to longer words $(cb)^{3}, a(bc)^{3}a, (bc)^3$.

\subsubsection{A presentation for space group $I\bar42d$ (No. 122)}
Consider the following generating set $\{ a, b, c \}$ for the $I\bar42d$ space group \cite{Fischer93}:

\begin{center}
\begin{tabular}{l}
$a=2 (x, 1/4, 1/8): x, 1/2-y, 1/4-z$; \\
$b=2(1/4, y, -1/8): 1/2-x, y, -1/4-z$; \\
$c=\bar{4}_{z}(000): y,-x,-z$.
\end{tabular}
\end{center}

\noindent The $I$-lattice is generated by $(ab)^2, abc^2, (ac)^2$. The presentation of the point group $\bar42m$ on the linear parts of $\{a, b, c\}$ is as follows:

\begin{center}
$\langle {\bar{a}}^2, {\bar{b}}^2, \bar{b}{\bar{c}}^{2}\bar{a}, (\bar{a}{\bar{c}}^{-1})^2, {\bar{c}}^2\bar{b}\bar{a} \rangle$.
\end{center}

\noindent The last three relators give rise to non-trivial translations once the point-group operations are lifted to their space-group counterparts:

\begin{center}
$bc^{2}a * (ac)^2 = (ac^{-1})^2 * (abc^{2})^{-1} = c^{2}ba * abc^{2}$,
\end{center}

\noindent where asterisks are merely a guide for the eyes to separate the original relators of the point group from the associated translations. As a result, the relators $bc^{3}ac$, $(ac^{-1})^2c^{-2}ba$, $c^{4}$ are found. Keeping in mind that $c^{-1} = c^{3}$ and clearly $c = c^{-3}$ we obtain $bc^{-1}ac, ac^{-1}acba, c^{4}$.  The second relator is identical to the first one since $ac^{-1}acba = c^{-1}acbaa = c^{-1}acb = bc^{-1}ac$.

Now we have to examine the conjugation action of each generator $a, b, c$ on the basis vectors of the $I$-lattice (to simplify the expressions, we introduce the notation: $t_{x}=(ab)^{2}, t_{y} =  abc^{2}, t_{z} = (ac)^{2}$):

\begin{center}
\begin{tabular}{l}
$at_{x}a = t_{x}^{-1}; bt_{x}b = t_{x}^{-1}; c^{-1}t_{x}c = t_{x}^{-1}$;\\

$at_{y}a = t_{z}^{-1}; bt_{y}b = t_{x}^{-1}t_{z}; c^{-1}t_{y}c = t_{z}^{-1}$;\\

$at_{z}a = t_{y}^{-1}; bt_{z}b = t_{x}^{-1}t_{y}; c^{-1}t_{z}c = t_{x}^{-1}t_{y}$.
\end{tabular}
\end{center}

\noindent A new relator is derived from the relation $c^{-1}t_{x}c = t_{x}^{-1}$. We have: 

\begin{center}
 $c^{-1}(abab)c(abab) = abcaba(bc^{-1}a)b$.
\end{center}

\noindent In the expression on the right there is a combination of letters in brackets that resembles the relator $bc^{-1}ac$. Now we insert a trivial word $cc^{-1}$ and achieve an elegant cancellation: $abcaba(bc^{-1}ac)c^{-1}b = abcabac^{-1}b$.

It can be verified that the commutators of the generators of the translation lattice need not be considered, and we finally obtain the following \emph{short} presentation for the $I\bar42d$ space group:

\begin{center}
$\langle a, b, c | a^{2}, b^{2}, c^{4}, bc^{-1}ac, abcabac^{-1}b \rangle$.
\end{center}

\begin{center}
\begin{figure}[ht]
\centering
\includegraphics{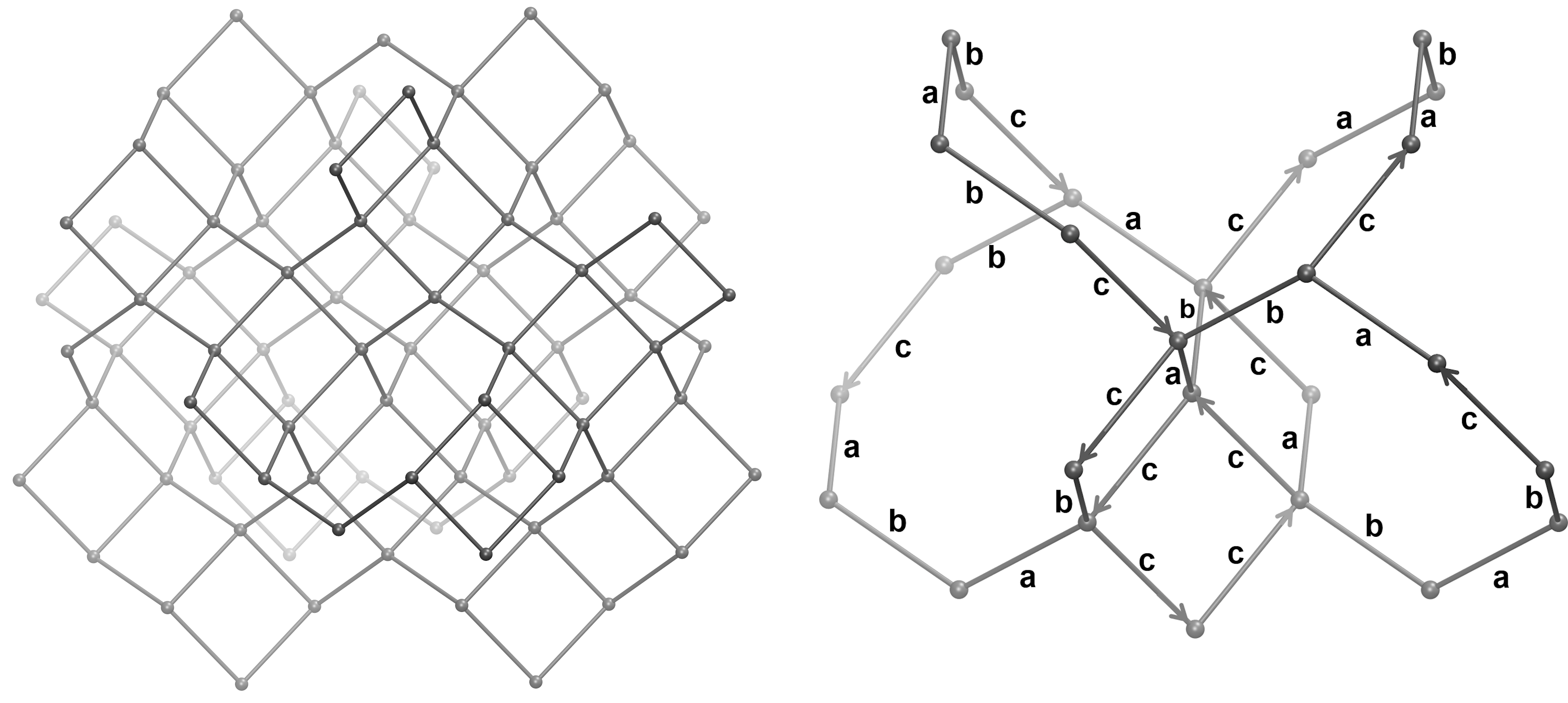}
\caption{Left: a finite portion of the GIS network. Right: cycles corresponding to relators (arrows distinguish between $c$~and~$c^{-1}$).}
\label{f:gispres}
\vspace{-2em}
\end{figure}
\end{center}

The Cayley graph for this generating set is isomorphic to the bond network in the zeolite gismondine (GIS) (Figure~\ref{f:gispres}). The number of symmetry-equivalent cycles sharing a vertex can be calculated from a simple formula $c \times N_{c}/N_{v}$ where $c$ is the cycle length, $N_{c}$ and $N_{v}$ stand for the number of these cycles and total number of vertices (of the whole graph) per unit cell, respectively.\footnote{We borrowed this formula from Fischer and Koch (1996) \cite[p.~2134]{FK96} who used it for counting the number of edge-sharing polygons in a tiling on a minimal surface.} In this way we obtain (orbit cardinalities refer to the conventional $I$-cell):

\begin{center}
\begin{tabular}{l}
$4 \times (4/16) = 1$ \emph{square} cycle with symmetry $\bar4$ (relator $c^{4}$);\\
$4 \times (8/16) = 2$ 4-cycles with 2-fold symmetry (relator $bc^{-1}ac$);\\
$8 \times (8/16) = 4$ 8-cycles with 2-fold symmetry (relator $abcabac^{-1}b$).
\end{tabular}
\end{center}


\newpage

\section{Results and discussion}

The proposed approach has been implemented in the $\GAP$ programming language (making use of the packages \emph{Cryst} \cite{Cryst2019} and \emph{Polycyclic}\cite{Polyc2020}) and applied to many generating sets that are of interest in crystallography.

For checking purposes we computed regular subgroups for the primitive cubic lattice ({\bf pcu}) along with their presentations. The resulting list contains 298 subgroups (up to conjugacy in $Aut({\bf pcu})=Pm\bar3m$) and is complete as the comparison with previously published lists shows \cite{Prok91, Delgado99, Kost2021}. It is \emph{a priori} clear that relators in \emph{shortest} presentations of these groups are of lengths 2 and 4 (the latter correspond to square cycles), and this holds indeed for more than $80\%$ of cases. The remaining presentations contain also longer relators (however, mostly not longer than 8) mainly because translations correspond to relatively long words and therefore, Tietze transformations are not always successful in deriving shortest possible relators. We anyway consider this result as quite satisfactory, and have not pursued a better performance. Moreover, the results of Tietze transformations depend on some random choices made in the substring searches and substitutions, -- for example, permuting generators and/or choosing another basis for translation group help to obtain shortest relators in many cases.

Since the amount of computed material is rather large and mainly suitable for supplementary information, we include here the results that seem bo be of special interest. Table~\ref{t:table1} lists presentations of regular subgroups for a number of high-symmetry 2- and 3-periodic graphs. For each group generators $a, b, c, ...$ should be interpreted as abstract generators, and for different groups they should not be mixed up with each other (this applies to all our tables). For 2- and 3-periodic graphs we use three-letter (lower-case bold) codes from the RCSR database of M.~O'Keeffe (https://rcsr.anu.edu.au), the only exception being the $2$-periodic minimal graph 2(3)2 for which the original name \cite{Klee92} is used.\footnote{We prefer here the term \emph{minimal graph} instead of the original \emph{minimal net} \cite{Klee92} to be consistent with the group- and graph-theoretic language of this paper.} For each graph its maximum-symmetry space group (in the international notation) is given that is isomorphic to the $Aut(\Gamma)$, oriented site-symmetry symbols \cite{IT} give information on vertex stabilizers. For the honeycomb tessellation ({\bf hcb})\footnote{To be more precise, the 1-skeleton of the honeycomb tessellation is meant here.} of the Euclidean plane we prefer however to use a slightly larger automorphism group that includes also a reflection in the plane that stabilizes it pointwise. For the 2(3)2 graph (that is non-planar) its $Aut(\Gamma)$ is isomorphic to the layer group $p\bar4m2$ (that is in turn isomorphic to the wallpaper group $p4mm$). To characterize local configuration of cycles around a vertex in a graph, we use the symbols -- that closely resemble those of Schl\"afli -- which give information on all \emph{strong rings} shared by the vertex, with superscripts referring to the number of rings of the same size. For example, in this notation the symbol for the primitive cubic lattice is $4^{12}$, for the diamond structure ({\bf dia}) -- $6^{12}$. This notation is relatively short, especially if hundreds or thousands of strong rings are to be specified.

Table~\ref{t:table2} contains complete lists of regular groups for twenty-two vertex-transitive zeolite structures (where they are designated by international three-letter upper-case codes). Table~\ref{t:table3} contains a presentation of one (arbitrarily chosen) group per a structure.

As an application to higher dimensions, we determined automorphism groups for vertex-transitive 4-periodic minimal graphs using the method of barycentric placement \cite{Delgado2004}. For each minimal graph Table~\ref{t:table4} contains its name in the catalogue \cite{Klee92} along with the name(s) of its 3-periodic quotients (which bring in useful associations from a crystallographic context), local configuration of strong rings around a vertex (which were computed partly following \cite{Eon2016}), the number of a four-dimensional space group in \cite{BBNWZ}, the structure of vertex stabilizer and transitivity $pqr$. The latter parameter indicates the number of vertex ($p)$, edge ($q$) and strong ring ($r$) orbits in $Aut(\Gamma)$. Presentations of regular groups for tri- and tetravalent 4-periodic minimal graphs are collected in Table~\ref{t:table5}.

\newpage

{\tiny
\begin{table}[H]
\begin{center}
\caption{Presentations of regular groups for some 2- and 3-periodic graphs}
\label{t:table1}
\begin{tabular}{|l|}
\hline
\multicolumn{1}{|c|}{$6^{3}\text{-{\bf hcb}} \, p6/mmm - \, \bar6m2$} \\ 
$p\bar1 \cong p112 \cong p112/b \text{ (twice)} = \langle a, b, c |  a^2, b^2, c^2, (abc)^2 \rangle$ \\
$cm11 \cong c121 \cong p2_{1}mn \cong p2an = \langle a, b |  b^2, ba^2ba^{-2} \rangle$ \\
$p12_{1}/a1 \text{ (twice)} \cong pba2 \cong p2_{1}2_{1}2 = \langle a, b |  b^2, (ba^2)^2 \rangle$ \\
$p12_{1}/m1 \cong p12/a1 \cong p2_{1}22 \cong pma2 \cong pbaa \cong pmab = \langle a, b, c |  a^2, b^2, c^2, cbcaba \rangle$ \\
$p3m1 \cong p312 = \langle a, b, c |  a^2, b^2, c^2, (cb)^3, (ac)^3, (ba)^3 \rangle$ \\
$p6 \cong p\bar3 = \langle a, b |  a^2, b^6, (ab)^3 \rangle$ \\ 
\hline
\multicolumn{1}{|c|}{$8^{4}\text{-{\bf 2(3)2}} \, p\bar4m2 - \, 2mm.$} \\ 
$c121 = \langle a, b |  a^2, (b^{-1}a)^2(ba)^2 \rangle; p2_{1}22 = \langle a, b, c |  a^2, b^2, c^2, (bc)^2(ac)^2 \rangle$ \\
$p\bar4b2 = \langle a, b |  b^2, (a^{-1}bab)^2, (ba)^4 \rangle; p\bar4m2 = \langle a, b, c |  a^2, b^2, c^2, (bcac)^2, (bc)^4, (ca)^4 \rangle$ \\ 
\hline
\multicolumn{1}{|c|}{$10^{15}\text{-{\bf srs}} \, I4_{1}32 - \, .32$} \\ 
$I2_{1}2_{1}2_{1} = \langle a, b, c |  a^2, b^2, c^2, cacb(ca)^2ba, ca(cb)^2abcb \rangle$ \\
$I4_{1} = \langle a, b |  b^2, (ba^{-1})^3ba^3 \rangle; P4_{1}2_{1}2 = \langle a, b |  b^2, baba^{-1}(ba)^2a^2 \rangle$ \\ 
$P4_{1}22 = \langle a, b, c |  a^2, b^2, c^2, baca(bc)^3, (ab)^2cacbab \rangle$ \\ 
\hline
\multicolumn{1}{|c|}{$10^{10}\text{-{\bf ths}} \, I4_{1}/amd - \, 2mm.$} \\ 
$C2/c = \langle a, b, c |  a^2, b^2, c^2, (babca)^2, (cabca)^2  \rangle$ \\
$Fdd2 = \langle a, b |  a^2, b^{-1}ab(ba)^2b^{-2}a  \rangle$ \\
$P4_122 = \langle a, b, c |  a^2, b^2, c^2, babcabacba, cab(ca)^2cba  \rangle$ \\
$P4_12_12 = \langle a, b |  a^2, b^{-1}ab(ba)^2b^2a  \rangle$ \\
\hline
\multicolumn{1}{|c|}{$6^{12}\text{-{\bf dia}} \, Fd\bar3m - \, \bar43m$} \\ 
$P\bar1 = \langle a, b, c, d |  a^2, b^2, c^2, d^2, (bac)^2, (dab)^2, (cad)^2 \rangle$ \\
$C2 = \langle a, b, c |  b^2, c^2, bacba^{-1}c, ca^2ca^{-2}, ba^2ba^{-2} \rangle$ \\
$Cc = \langle a, b |  b^{-1}a^2ba^{-2}, b^{-1}ab^2a^{-1}b^{-1} \rangle$ \\
$P2/c = \langle a, b, c, d |  a^2, b^2, c^2, d^2, (bcd)^2, cabcba, dabdba, (acd)^2 \rangle$ \\
$P2_{1}/c = \langle a, b, c |  b^2, c^2, a^{-1}bcabc, (ca^2)^2, (aba)^2 \rangle$ \\
$P2_{1}/c = \langle a, b, c |  b^2, c^2, a^{-1}bcacb, (ca^2)^2, (aba)^2 \rangle$ \\
$C2/c = \langle a, b, c |  b^2, c^2, (ca^2)^2, ba^2ba^{-2}, (cab)^2 \rangle$ \\
$C2/c = \langle a, b, c, d |  a^2, b^2, c^2, d^2, (bac)^2, dabdba, dacdbc \rangle$ \\
$P222_{1} = \langle a, b, c, d |  a^2, b^2, c^2, d^2, bcdbdc, cabcba, dabdba, acdadc \rangle$ \\
$P2_{1}2_{1}2_{1} = \langle a, b |  b^{-1}a^2ba^2, a^{-1}b^2ab^2 \rangle; Pna2_{1} = \langle a, b |  b^{-1}a^2ba^{-2}, a^{-1}b^2ab^2 \rangle$ \\
$Pnc2 = \langle a, b, c |  b^2, c^2, a^{-1}bcacb, ca^2ca^{-2}, aba^{-2}ba \rangle$ \\
$Pnna = \langle a, b, c, d |  a^2, b^2, c^2, d^2, cdadcb, acdbdc, (dab)^2, cabcba \rangle$ \\
$Pbcn = \langle a, b, c |  b^2, c^2, cacbab, ba^2ba^{-2}, (ca^2)^2 \rangle$ \\
$P4_{1} = \langle a, b |  a^{-2}baba^{-1}, b^{-1}a^{-1}b^3a^{-1} \rangle$ \\
$P4_{1}2_{1}2 = \langle a, b, c |  a^2, c^2, cbcb^{-1}ca, b^{-1}abaca, b^3abc \rangle$ \\
$P3_{1}21 = \langle a, b, c |  a^2, b^2, c^{-1}bc^2bc^{-1}, abc^{-1}ac^{-2}, abcac^{-1}b \rangle$ \\ 
\hline
\multicolumn{1}{|c|}{$6^{8}\text{-{\bf nbo}} \, Im\bar3m - \, 4/mmm$} \\ 
$R3 = \langle a, b |  b^3a^3, ba^{-1}b^{-1}a^2b^{-1}a^{-1}b  \rangle \cong \langle a, b |  b^3a^3, (ba)^{3} \rangle$ \\
$R\bar3 = \langle a, b |  (b^2a)^2, (ba^2)^2, a^6, b^{-1}a^2b^{-1}a^{-1}b^2a^{-1}  \rangle \cong \langle a, b | a^6, b^6, (b^2a)^2, (ba^2)^2 \rangle$ \\
$R32 = \langle a, b, c, d |  a^2, b^2, c^2, d^2, bdadbc, dacbca, (ca)^3, bcdbadca  \rangle \cong$ \\
$\cong \langle a, b, c, d |  a^2, b^2, c^2, d^2, bdadbc, dacbca, (ca)^3, (bd)^3  \rangle$ \\
$R3c = \langle a, b |  baba^{-1}b^{-1}a^{-1}, a^{-1}b^2a^2b^{-2}a^{-1}  \rangle \cong \langle a, b |  baba^{-1}b^{-1}a^{-1}, (a^{-1}b)^3  \rangle$ \\
$R\bar3c = \langle a, b, c |  a^2, c^2, (cb^{-1}a)^2, (ca)^3, b^{-1}ab^2cb^{-1}, cb(ba)^2cb  \rangle \cong$ \\
$    \cong \langle a, b, c |  a^2, b^6, c^2, (cb^{-1}a)^2, (ca)^3, b^{-1}ab^2cb^{-1} \rangle$ \\ 
\hline
\multicolumn{1}{|c|}{$6^{6}.8^{40}\text{-{\bf qtz}} \, P6_{2}22 - \, 222$} \\ 
$P3_{2} = \langle a, b |  b^2a^{-3}b, a^{-2}(ba)^2b^{-2}  \rangle$ \\
$P6_{1} = \langle a, b |  a^{-1}b^{-1}abab^{-1}, b^{-1}a^{-1}b^3a^{-1}b^{-1}a  \rangle$ \\
$P6_{5}22 = \langle a, b, c |  b^2, c^2, a^{-2}ca^2b, (abc)^2, a^{-1}(a^{-1}c)^2bac  \rangle$ \\
$P3_{1}21 = \langle a, b |  (a^{-2}b)^2, (a^{-1}b^2)^2, b^{-1}a^2ba^2b^{-2}  \rangle$ \\
$P3_{1}12 = \langle a, b, c, d |  a^2, b^2, c^2, d^2, bacabd, acdbdc, dbadabdc  \rangle$ \\ \hline

\end{tabular}
\end{center}
\end{table}
}

{\tiny
\begin{table}[t]
\begin{center}
\caption{Regular groups for vertex-transitive zeolites}
\label{t:table2}
\begin{tabular}{|l|l|}
\hline
Zeolite & Groups \\ \hline
ABW & $C2/m, C2/c, I2_{1}2_{1}2_{1}, Ima2, Pnna \, \text{(twice)}, Pmna, Pnma$ \\ \hline
ACO & $Immm, I4/m, I422, I\bar4m2, P4/nnc, P4/mnc, P4_{2}/mmc, P4_{2}/nmc$ \\ \hline
AFI & $P6/mcc$ \\ \hline
ANA & $Ia\bar3, I\bar43d$ \\ \hline
ATN & $I422, I4mm, I\bar4m2, I\bar42m, P4/nnc, P4/nmm, P4_{2}/nnm, P4_{2}/nmc$ \\ \hline
ATO & $R\bar3m$ \\ \hline
BCT & \parbox[t]{8cm}{$C2/m, I222, Imm2, Pnnn, Pnnm, Pmmn, I4, I\bar4, P4/n$,
	$P4_{2}/n, P42_{1}2, P4_{2}22, P4nc, P4_{2}mc, P\bar42c, P\bar42_{1}c, P\bar4m2$,
	$P\bar4n2, P4/ncc, P4_{2}/nnm, P4_{2}/nmc, P4_{2}/ncm, I4_{1}/acd \, \text{(twice)}$} \\ \hline
BSV & $Ia\bar3d$ \\ \hline
CAN & $P\bar31c, P\bar3m1, P6_{3}22, P6_{3}mc$ \\ \hline
CHA & $R\bar3m$ \\ \hline
DFT & $Cccm, P4_{2}/m, P4_{2}22, P\bar42c, P4_{2}/mcm, P4_{2}/nnm, P4_{2}/mnm, P4_{2}/ncm$ \\ \hline
FAU & $Fd\bar3m$ \\ \hline
GIS & $Imma, I4_{1}/a, I4_{1}md, I\bar42d$ \\ \hline
GME & $P6_{3}/mmc$ \\ \hline
KFI & $Im\bar3m$ \\ \hline
LTA & $P432, P\bar43m, Fm\bar3m, Fm\bar3c$ \\ \hline
MER & $I4/mmm$ \\ \hline
MON & $Fddd, I4_{1}/a, I4_{1}22, I\bar42d$ \\ \hline
NPO & $P\bar3, P321, P31c, P\bar3c1, P6_{3}$ \\ \hline
RHO & $Im\bar3, I\bar43m, Pm\bar3m, Pn\bar3n$ \\ \hline
RWY & $I432, I\bar43m, Pn\bar3n, Pn\bar3m$ \\ \hline
SOD & $R\bar3, R3m, R\bar3m, R\bar3c, Fd\bar3, F432, F\bar43m, F\bar43c$ \\ \hline
\end{tabular}
\end{center}
\end{table}
}

{\tiny
\begin{table}[H]
\begin{center}
\caption{Presentations of regular groups for vertex-transitive zeolites}
\label{t:table3}
\begin{tabular}{|l|l|l|}
\hline
Zeolite & Schl\"afli & Presentation \\ \hline
ABW & $4^2.6^3.8^4$  & $Pnma = \langle a, b, c |  a^2, c^2, abab^{-1}, (cb^2)^2, cbcacb^{-1}ca  \rangle$ \\ \hline
ACO & $4^3.8^6$  & $P4/nnc = \langle a, b, c |  a^2, b^2, (ac)^2, c^4, (abcb)^2, (c^{-1}b)^4  \rangle$ \\ \hline
AFI & $4.6^{13}.12$  & $P6/mcc = \langle a, b, c, d |  a^2, b^2, c^2, d^2, (ca)^2, cdcada, abdbad, (bc)^3, (ba)^6  \rangle$ \\ \hline
ANA & $4^2.6^2.8^{16}$  & $I\bar43d = \langle a, b |  a^4, b^4, (ba)^3, ab^{-1}a^{-2}ba^{-1}b^{-2}  \rangle$ \\ \hline
ATN & $4^2.6^3.8^7$  & $P4_{2}/nmc = \langle a, b, c |  a^2, b^2, (bc)^2, ac^2ac^{-2}, (c^{-1}a)^2(ca)^2, (ba)^4  \rangle$ \\ \hline
ATO & $4.6^9.12^{20}$  & $R\bar3m = \langle a, b, c, d |  c^2, a^2, b^2, d^2, (dc)^2, (bc)^2ac, dabdba, (bd)^6  \rangle$ \\ \hline
BCT & $4.6^6.8^{12}$  & $I4_{1}/acd = \langle a, b, c |  a^2, b^2, c^4, (abc)^2, (bc^{-1}a)^2, c^{-1}b(ca)^2cb  \rangle$ \\ \hline
BSV & $4^3.6.12^{294}$  & $Ia\bar3d = \langle a, b |  (ba)^2, b^4, a^6, (a^{-1}b)^2ab^{-1}a^{-1}ba^{-2}b^{-2}  \rangle$ \\ \hline
CAN & $4^2.6^4.12^{20}$  & $P6_{3}mc = \langle a, b, c |  b^2, c^2, baba^{-1}, (bc)^3, ca^2ca^{-2}, (ca)^3(ca^{-1})^3  \rangle$ \\ \hline
CHA & $4^3.6.8^2$  & $R\bar3m = \langle a, b, c, d |  a^2, b^2, c^2, d^2, dbcb, (ac)^2, (cd)^3, (abad)^2  \rangle$ \\ \hline
DFT & $4^2.6^6.8^{10}$  & $P4_{2}/ncm = \langle a, b, c |  a^2, b^2, c^4, (ab)^2, (cbc)^2, cbcac^{-1}bc^{-1}a, (c^{-1}aca)^2  \rangle$ \\ \hline
FAU & $4^3.6^2.12$  & $Fd\bar3m = \langle a, b, c, d |  (ca)^2, dadb, (bc)^3, (ab)^3, (cd)^6  \rangle$ \\ \hline
GIS & $4^3.8^4$  & $I4_{1}/a = \langle a, b |  (ab)^2, b^4, (b^{-1}a^3)^2  \rangle$ \\ \hline
GME & $4^3.6.8^2.12$  & $P6_{3}/mmc = \langle a, b, c, d |  a^2, b^2, c^2, d^2, (db)^2, (cb)^2, (ca)^2, (cd)^3, (abad)^2, (da)^6  \rangle$ \\ \hline
KFI & $4^3.6.8^2$  & $Im\bar3m = \langle a, b, c, d |  a^2, b^2, c^2, d^2, cada, (bd)^2, (cd)^3, (cb)^4, (ab)^4  \rangle$ \\ \hline
LTA & $4^3.6^2.8$  & $Fm\bar3c = \langle a, b, c |  a^2, b^2, c^4, acac^{-1}, (bc^{-1})^3, (ba)^4  \rangle$ \\ \hline
MER & $4^3.8^4$  & $I4/mmm = \langle a, b, c, d |  b^2, d^2, a^2, c^2, (cb)^2, (ad)^2, (db)^2, (acab)^2, (cd)^4, (ca)^4  \rangle$ \\ \hline
MON & $4.5^5.8^6$  & $I\bar42d = \langle a, b, c |  a^2, b^2, c^4, babc^{-2}, c^{-1}ac^{-1}bcacb, (c^{-1}a)^4  \rangle$ \\ \hline
NPO & $3.6^6.12^{40}$  & $P\bar3c1 = \langle a, b, c |  a^2, b^2, c^3, bcbaca, c^{-1}aca(c^{-1}b)^4  \rangle$ \\ \hline
RHO & $4^3.6.8^2$  & $Pn\bar3n = \langle a, b |  (a^{-1}b^{-1})^2, a^4, b^6, (a^{-1}b)^4  \rangle$ \\ \hline
RWY & $3^3.8.12^2$  & $Pn\bar3m = \langle a, b, c, d |  a^2, b^2, c^2, d^2,  dca, (ba)^4, (bc)^6, (bd)^6 \rangle$ \\ \hline
SOD & $4^2.6^4$  & $F\bar43c = \langle a, b |  a^4, b^4, (ab)^3, (ab^{-1})^3  \rangle$ \\ \hline
\end{tabular}
\end{center}
\end{table}
}

{\tiny
\begin{table}[t]
\begin{center}
\caption{Symmetry properties of vertex-transitive 4-periodic minimal graphs}
\small Notation for four-dimensional space groups as in \cite{BBNWZ}.
\label{t:table4}
\begin{tabular}{|c|c|c|c|c|c}
\hline
Name & Schl\"afli & Transitivity & $Stab_{Aut(\Gamma)}(x)$ & $Aut(\Gamma)$ \\ \hline
1(8)1 & $4^{24}$ & 111 & $C_{2}^{4} \rtimes S_4$ & 32/21/1/1 \\ \hline
2(5)1 & $6^{30}$ & 111 & $S_5$ & 31/7/2/1 \\ \hline
2(5)2 & $6^{21}.8^{60}$ & 124 & $C_2 \times C_2 \times S_3$  & 20/20/2/1 \\ \hline
2(5)3 & $4^4.8^{16}$ & 122 & $D_4 \times D_4$ & 32/17/1/1 \\  \hline
3(4)1 & $8^{40}$ & 112 & $D_4 \times C_2$ & 25/8/3/3\\ \hline
3(4)3 & $8^{16}$ & 112 & $D_4 \times C_2$ & 25/8/1/3 \\  \hline
6(3)1 & $12^{12}.14^{42}$ & 112 & $D_3 \times C_2$ & 29/8/2/2 \\ \hline
6(3)2 & $12^{12}.14^{21}$ & 123 & $C_2$ & 14/7/2/2 \\ \hline
6(3)4 & $12^{6}.14^{28}$ & 122 & $C_2 \times C_2 \times C_2$ & 25/8/3/3 \\
\hline
\end{tabular}
\end{center}
\end{table}
}

{\tiny
\begin{table}[H]
\begin{center}
\caption{Presentations of regular groups for some 4-periodic graphs}
\label{t:table5}
\begin{tabular}{|l|}
\hline
\multicolumn{1}{|c|}{$8^{40}$ -- {\bf 3(4)1 `qtz'}} \\ 
$(8/1/1/2) = \langle a, b |  b^{-3}a^{-1}b^3a, aba^{-2}bab^{-2}  \rangle$ \\
$(14/1/1/2) = \langle a, b |  b^{-1}a^{-2}b^{-1}a^4, bab^{-2}a^{-1}b^{-1}a^2  \rangle$ \\
$(14/7/1/6) = \langle a, b, c |  b^2, c^2, baca^{-1}(a^{-1}c)^2, abca^{-1}ca^2c  \rangle$ \\
$(14/7/1/6) = \langle a, b, c |  b^2, c^2, a^2ca^4b, ca^{-2}ca^{-1}bca  \rangle$ \\
$(14/3/1/4) = \langle a, b |  ba^2b^{-1}a^{-1}b^2a, ba^{-2}b^{-1}ab^{-2}a^{-1}  \rangle$ \\
$(8/4/1/2) = \langle a, b, c, d |  a^2, b^2, c^2, d^2, cadabdad, bc(db)^2ad  \rangle$ \\
$(8/4/1/2) = \langle a, b, c |  b^2, c^2, (a^{-3}c)^2, bca^{-2}ca^{-1}ba, ca^{-1}ba^{-2}bca  \rangle$ \\
\hline
\multicolumn{1}{|c|}{$12^{6}.14^{28}$ -- {\bf 6(3)4 `bto'}}\\ 
$(14/3/3/2) = \langle a, b |  b^2, aba^2(ba^{-2})^2ba, (ba)^2aba^{-1}baba^{-2}ba^{-1}  \rangle$ \\
$(8/4/1/2) = \langle a, b, c |  a^2, b^2, c^2, bcab(cba)^2ca, (b(ac)^3)^2  \rangle$ \\
$(14/7/1/2) = \langle a, b, c |  a^2, b^2, c^2, b(cba)^2cabca, (ca)^2b(ca)^2cbaca  \rangle$ \\
$(14/7/1/6) = \langle a, b |  b^2, aba^2(ba^{-2})^2ba, ((ba)^2aba^{-1})^2  \rangle$ \\
\hline
\multicolumn{1}{|c|}{$12^{12}.14^{42}$ -- {\bf 6(3)1 `twt/pcu-h'}} \\ 
$(8/4/1/2) = \langle a, b, c |  a^2, b^2, c^2, (ca)^2(bc)^2(ab)^2, (bcbacbc)^2  \rangle$ \\
$(14/1/2/2) = \langle a, b |  b^2, ba^{-4}b(aba)^2, a^{-1}ba^{-2}(ba)^3aba^{-1}b  \rangle$ \\
$(14/7/1/7) = \langle a, b, c |  a^2, b^2, c^2, c(acb)^2c(ba)^2, (ab)^2acabcacbac, cbcabc(ba)^2bcba  \rangle$ \\
$(14/7/1/6) = \langle a, b |  b^2, a^{-1}ba^{-3}(a^{-1}b)^2a^2b, ba^{-1}(ba^2)^2a^3ba^2  \rangle$ \\
\hline
\multicolumn{1}{|c|}{$12^{12}.14^{21}$ -- {\bf 6(3)2 `eta/etb'}} \\ 
$(14/1/2/2) = \langle a, b |  b^2, ba^3(ba^{-1})^2a^{-2}ba, ba^2(ba^{-1})^4ba^2  \rangle$ \\
$(14/3/3/2) = \langle a, b |  b^2, (ba^{-1})^2a^{-2}(ba)^2a^2, b(a^{-1}ba^{-1})^2b(aba)^2  \rangle$ \\
$(14/7/1/7) = \langle a, b, c |  a^2, b^2, c^2, a(bc)^2ac(ba)^2c, cbac(bcba)^2ba  \rangle$ \\
$(14/7/1/2) = \langle a, b, c |  a^2, b^2, c^2, a(cb)^2ac(ab)^2c, bac(ba)^4bca  \rangle$ \\
\hline
\end{tabular}
\end{center}
\end{table}
}

\subsection{2-periodic graphs}

The presentations for regular subgroups of the {\bf hcb} tessellation provide an elegant way to clearly see the isomorphism of certain layer groups regarded as abstract groups. Altogether in this case there are twenty-two conjugacy classes of regular subgroups (two `isosymbolic' subgroups of types $p112/b$ and $p12_1/a1$ occur twice each as indicated in Table~\ref{t:table1}) which belong to six (abstract) isomorphism classes that fully agrees with the results of Gr\"unbaum and Shephard (1977) \cite{Gruen77}. A comparison with the subgroup data of Koch and Fischer (1978) \cite{KF78} shows perfect agreement apart from a minor discrepancy: in their table on p.146 the symbol $p12/m(1)$ should be changed to $p12_{1}/m(1)$ (that this is only a misprint can be seen from the given symmetry operations which indeed generate the correct $p12_{1}/m(1)$ group and not another one).\footnote{In this remark we use the nomenclature for subperiodic groups as in Koch and Fischer (1978).}

The presentations for regular subgroups of the 2(3)2 graph are believed to be new. Incidentally, our presentation for the group $c121$ appears on p.44 in the book of Coxeter and Moser (1980) \cite{Coxeter80} in connection with the wallpaper group $cm$ (that is isomorphic to the layer group $c121$). However, the presentation is not discussed in any detail and is even not numbered.

Besides 2-periodic graphs which are familiar as ornamental tilings of the Euclidean plane, those which resemble locally 3-periodic graphs, and therefore give rise to the dilemma which structure should be preferred in experiment, are particularly interesting from a crystallographic perspective. Mitina and Blatov (2013) \cite{Blatov13} described one of the layered coordination polymers (CSD code HOWVAR) as based on the graph that is locally similar to the 3-periodic $10^{10}$-{\bf ths} graph, albeit it contains also 12-rings. The original crystal structure belongs to the space group $P2_{1}/c$, individual layer (layer group $p1c1$) contains four vertex orbits (which correspond to molecular barycenters). The automorphism group of the layer is isomorphic to the wallpaper group $c2mm$ that is vertex-transitive with trivial vertex stabilizers meaning that the layer is a Cayley graph for this group. A closer inspection has shown that this graph is indeed a quotient of the {\bf ths} graph by some rank-1 translation subgroup.

\begin{center}
\begin{figure}[ht]
\centering
\includegraphics{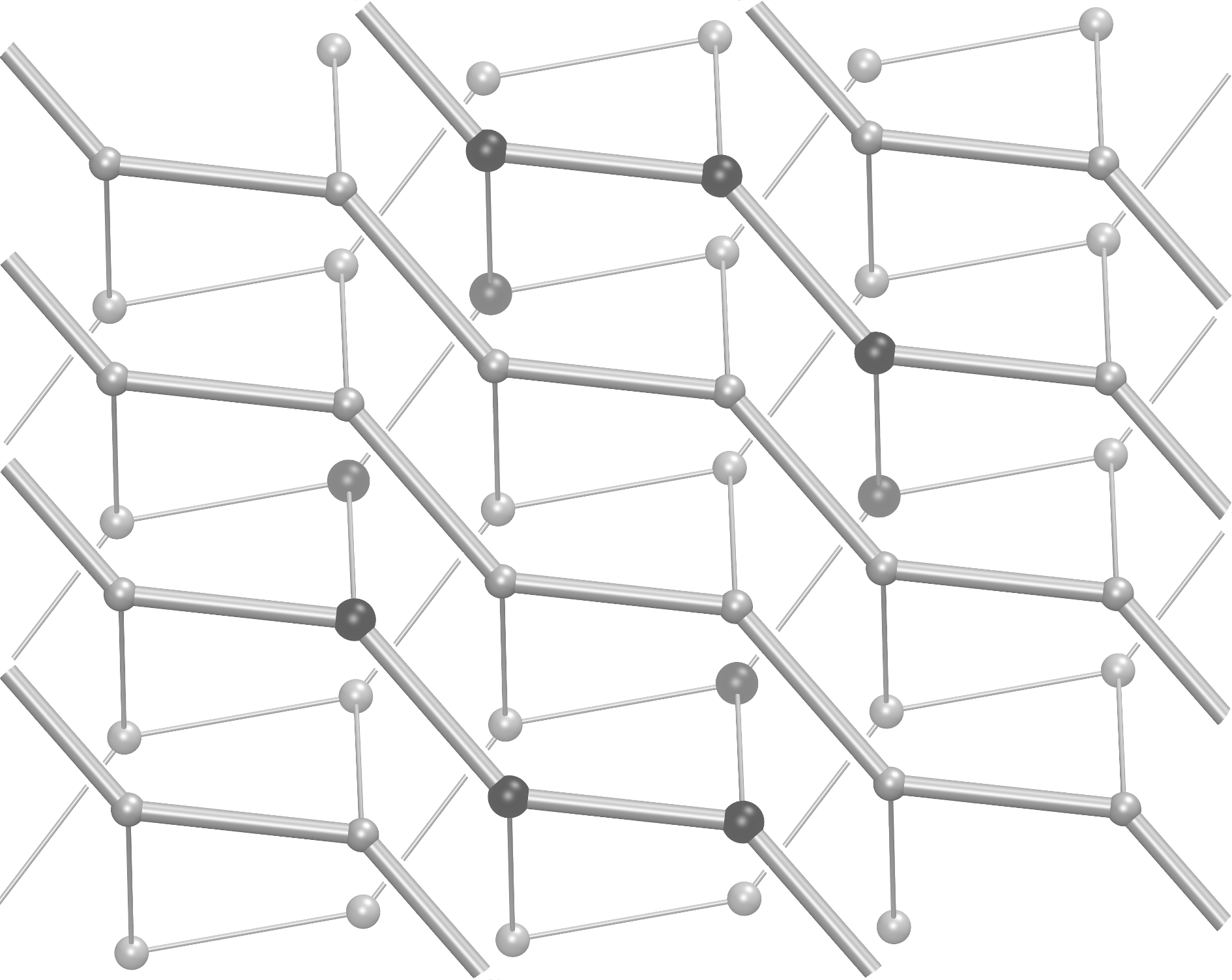}
\caption{The 2-periodic quotient ($10^{10}.12^6$)  of the {\bf ths} embedded with layer-group symmetry $p2_122$. The vertices of the 10-ring are marked.}
\label{f:thspmg}
\end{figure}
\end{center}

To explore 2-periodic quotients of the {\bf ths} in more detail, we enumerated the quotients by those translations which correspond to geodesics of length 12 using the projection method of Eon \cite{Eon2011}. To this end, only translations with coprime components (referred to a primitive basis) were considered because otherwise the quotients would possess bounded automorphisms of finite order that is undesirable for many reasons. As a result, we identified only three vertex-transitive representatives (Table~\ref{t:ths}), including the one observed in experiment (highlighted in bold). The graphs obtained are non-planar and allow for piecewise-linear embeddings in $\R^3$ with symmetry of layer groups that are isomorphic to their automorphism groups (Figure~\ref{f:thspmg}). The embeddings (see supporting information for coordinates) show self-catenation as expected. 

{\tiny
\begin{table}[h]
\caption{Vertex-transitive 2-periodic quotients of the {\bf ths} graph}
\label{t:ths}
\small \raggedright \hspace{2em} Translations refer to the conventional $tI$-cell.\\
\small \raggedright \hspace{2em} TD\textsubscript{10} is the number of vertices in a ball of radius 10.
\begin{center}
\begin{tabular}{|l|l|l|l|}
\hline
Translation & Schl\"afli & TD\textsubscript{10} & Presentation \\ \hline
${\bf [5/2, 5/2, 1/2]}$ & ${\bf 10^{10}.12^3}$  & 424 & $c2mm = \langle a, b, c | a^2, b^2, c^2, (abacb)^2, (cbacb)^2, (ab(ac)^2)^2 \rangle$ \\ \hline
$[2, 2, 1]$          &  $10^{10}.12^6$  & 445 & $p2mg = \langle a, b, c | a^2, b^2, c^2, ((ac)^2b)^2, (bacbc)^2, (abcb)^2cbcb \rangle$ \\  \hline
$[1/2, 1/2, -3/2]$   &  $10^{10}.12^9$  & 460 & $c2mm = \langle a, b, c | a^2, b^2, c^2, (cabcb)^2, ba(bc)^4ba \rangle$ \\  \hline
\end{tabular}
\end{center}
\end{table}
}

\noindent Wallpaper groups in Table~\ref{t:ths} also naturally arise as the quotient groups $N_{Aut({\bf ths})}(\langle t \rangle)/\langle t \rangle$ where $t$ is the respective translation \cite{Sene88}. The normalizer in all three cases is $C2/c$ (\emph{cf.} Table~\ref{t:table1}). The number of 12-rings at the vertex is equal to the number of geodesics of length $12$ in the original {\bf ths} which are factored out in the quotient. Note that the presentations contain only one relator of length 12, although the number of orbits of 12-rings can be larger (as follows from Schl\"afli symbols $10^{10}.12^6$ and $10^{10}.12^9$). We shall return to this point in the next section. 

\subsection{3-periodic graphs}
Presentations for regular subgroups of 3-periodic graphs from Tables 1 and 3 seem not to have appeared in the literature before, apart from some presentations of fixed-point-free groups in \cite{Molnar87} and for $P2_{1}2_{1}2_{1}$ in \cite{Ban19}. However, presentations associated with many tetravalent graphs should be implicit in the derivation of 4-regular vertex-transitive tilings of $\R^3$ undertaken by Delgado-Friedrichs and Huson \cite{Delgado2000}.

Having in mind possible applications in structural chemistry, we would like to emphasize two non-conjugate subgroups of type $P2_{1}/c$ of the {\bf dia} graph. That both subgroups are indeed non-conjugate in $Aut({\bf dia})=Fd\bar3m$ can be seen from the fact that one generating set contains a screw rotation as the element of infinite order, whereas another one contains a glide plane. This subtle point can also be read off from the tables by Sowa (2019) \cite{Sowa2019}.

\begin{center}
\begin{figure}[h]
\centering
\includegraphics{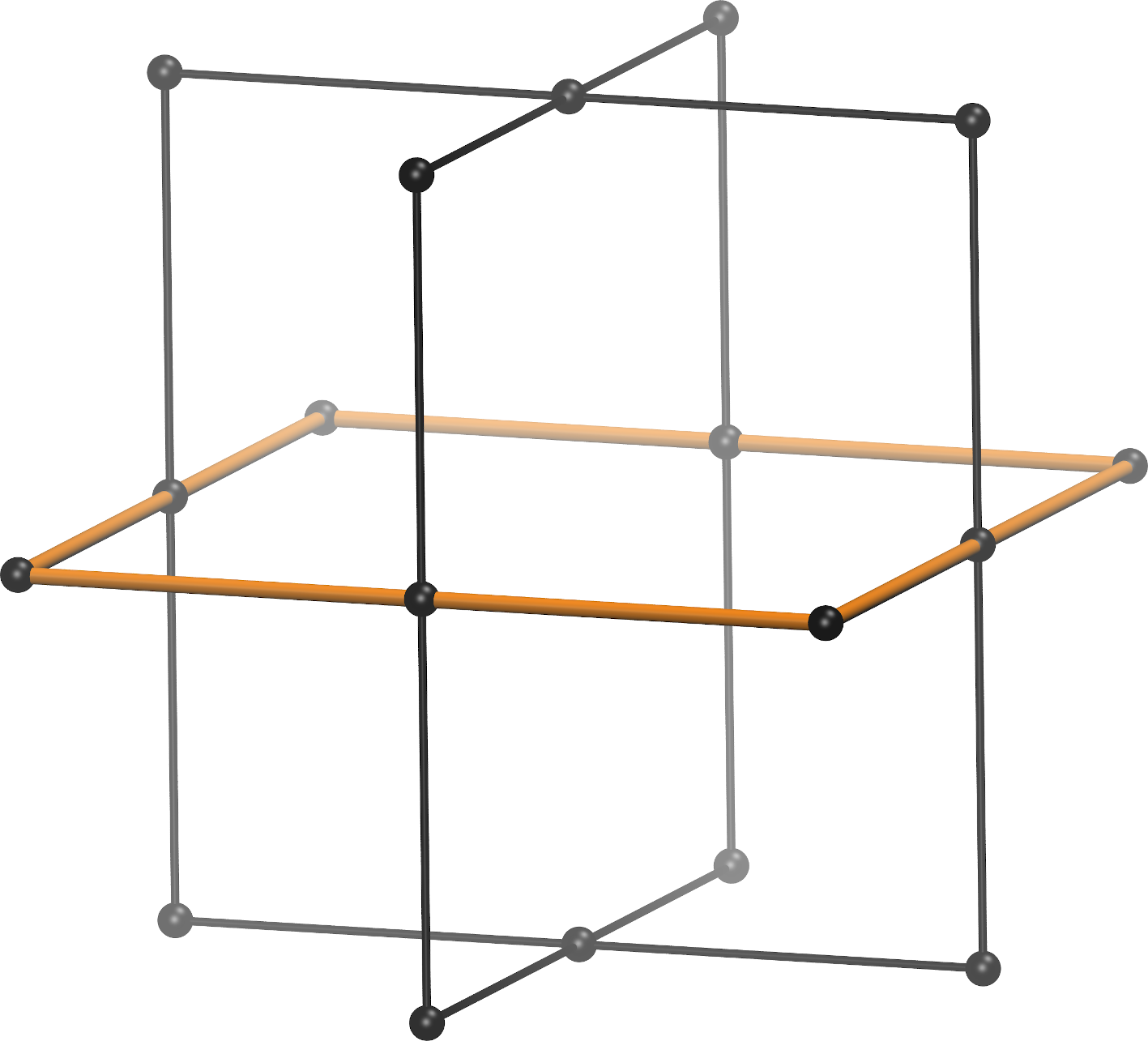}
\caption{Cycles in {\bf nbo}. The 8-cycle is emphasized.}
\label{f:nboring}
\vspace{-1em}
\end{figure}
\end{center}

Since strong rings in the high-symmetry graphs from Table~\ref{t:table1} and also for zeolites (Table~\ref{t:table3}) are well understood \cite{Blatov07}, they provide good examples for checking whether computed presentations are indeed the shortest. In all cases but one our program was indeed able to determine shortest possible presentations (sometimes permuting generators was necessary, \emph{e.g.} for FAU and RWY). The unlucky case is the {\bf nbo} 3-periodic graph for which our presentations contain relators of length~8. Figure~\ref{f:nboring} demonstrates that the corresponding 8-cycle can be decomposed into the sum of four shortest 6-cycles. Such instances cannot be always overcome by standard Tietze transformations routines.

Along with highly-symmetric graphs with a transparent `cycle structure' (as can be seen from short presentations), there are a number of vertex-transitive 3-periodic graphs graphs with hundreds and thousands of strong rings at the vertex. First examples of this kind were found among sphere-packing graphs (W. Fischer and collaborators). To take a closer look at such exceptional cases, we used the \emph{julia} package \emph{PeriodicGraphs} \cite{Zoub24} to compute strong rings for all published sphere-packing graphs (see supporting information). The examples of exceptional sphere-packing graphs (for which the original Fischer nomenclature \cite{Fischer71} is used) and presentations for their vertex-transitive groups are collected in Table~\ref{t:monsters}.

{\tiny
\begin{table}[h]
\begin{center}
\caption{Sphere-packing graphs with exceptionally high number of strong rings}
\label{t:monsters}
\begin{tabular}{|l|}
\hline
\multicolumn{1}{|c|}{${\bf 4/3/c22 - 3.4.6^{12}.20^{5120}}$} \\ 
$I4_{1}32 = \langle a, b, c | a^2,  b^3, c^2, (ac)^2, cbcaba, ab^{-1}abab^{-1}(ab^{-1}(ab)^2)^2ab^{-1} \rangle$ \\
\hline
\multicolumn{1}{|c|}{${\bf 4/3/c23 - 3.4^{2}.6^{6}.20^{415}}$} \\ 
$I4_{1}32 = \langle a, b, c |  a^2, c^2, b^3, (cb^{-1})^2, (abc)^2, (ab^{-1}ac)^2ab^{-1}(ab)^2ab^{-1}(ac)^2  \rangle$ \\
\hline
\multicolumn{1}{|c|}{$ {\bf 4/3/c24 - 3.4^{2}.6^{6}.20^{310}}$}\\ 
$I4_{1}32 = \langle a, b, c |  b^2, c^2, a^3, (a^{-1}b)^2, bcbaca^{-1}, cbcacbca^{-1}(cacbca)^2  \rangle$ \\
\hline
\multicolumn{1}{|c|}{$ {\bf 5/4/c6 - 4^{12}.12^{1458}}$} \\ 
$Pn\bar3n = \langle a, b, c |  c^2, (b^{-1}c)^2, ca^2b, b^4, ab(a^{-1}c)^5  \rangle$ \\
\hline
\multicolumn{1}{|c|}{$ {\bf 4/4/t60 - 4^{2}.6^{3}.8^{6}.14^{553}}$} \\ 
$I4_{1}/acd = \langle a, b, c | a^2, b^2, (ac)^2, (bc^{-2})^2, (c^{-1}b)^2(cb)^2, ac^{-1}(ba)^3c^{-2}bc^3 \rangle$ \\
\hline
\multicolumn{1}{|c|}{$ {\bf 4/4/h14 - 4^{2}.6^{6}.18^{3615}}$} \\ 
$R\bar3c = \langle a, b, c, d |  a^2, b^2, c^2, d^2, dbda, (ca)^2cb, (bdc)^6  \rangle$ \\
\hline
\multicolumn{1}{|c|}{$ {\bf 4/4/h17 - 4^{2}.6^{6}.18^{3213}}$} \\ 
$P6_{2}22 = \langle a, b, c, d | a^2, b^2, c^2, d^2, (da)^2, (ac)^2, (bd)^2bc, ba(cbab)^2adbabdba  \rangle$ \\
\hline
\multicolumn{1}{|c|}{$ {\bf 5/4/h16 - 4^{7}.8^{28}.16^{1820}}$} \\ 
$P6_{2}22 = \langle a, b, c, d |  a^2, b^2, c^2, (db)^2, (cb)^2, (cd)^2, (ad)^2, adbd^{-1}(ab)^2, d^{-1}(ac)^2d^{-3}(ac)^4 \rangle$ \\
\hline
\multicolumn{1}{|c|}{$ {\bf 5/4/h19 - 4^{6}.6^{20}.18^{274428}}$} \\ 
$P6cc = \langle a, b, c |  a^2, acac^{-1}, c^2b^{-2}, (bc^{-1})^3, b^{-1}c^{-1}abacbc^{-1}abacb^{-1}c^{-1}abc^{-1}a  \rangle$ \\
$P\bar31c = \langle a, b, c, d, f |  a^2, b^2, c^2, d^2, f^2, bdfc, fada, cfcdbf, badbafcfacabd(ba)^2d  \rangle$ \\
$P622 = \langle a, b, c, d, f |  a^2, b^2, c^2, d^2, f^2, bdfc, (fa)^2, (ad)^2, fcfbdc, babfacabdbacf(ab)^2d  \rangle$ \\
$P\bar3c1 = \langle a, b, c |  a^2, c^2b^{-2}, (ac)^2, b^{-1}c(bc^{-1})^2, b^{-1}c^{-1}(ab)^2cbab(c^{-1}ab^{-1})^2c^{-1}a  \rangle$ \\
\hline
\multicolumn{1}{|c|}{$ {\bf 5/4/h29 - 4^{6}.5^{5}.7^{7}.17^{6188}}$} \\ 
$P6_{1}22 = \langle a, b, c, d | a^2, b^2, c^2, (bd^{-1})^2, (cd)^2, (ad)^2, d(cb)^2, d(ba)^3, adcacd^{-3}(ac)^4d^{-1} \rangle$ \\
\hline
\end{tabular}
\end{center}
\end{table}
}

\begin{center}
\begin{figure}[H]
\centering
\includegraphics{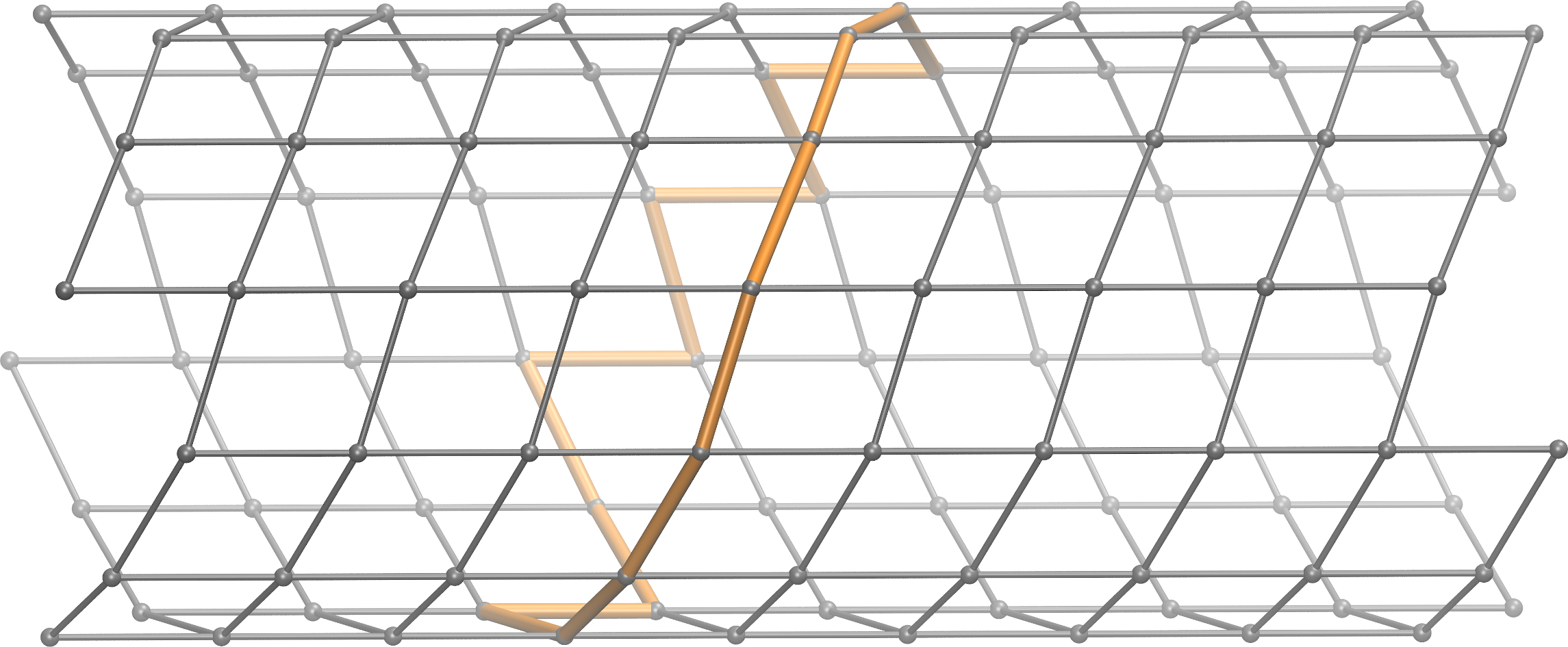}
\caption{Square lattice rolled along the vector (4, 12). One of the 16-rings is marked.}
\label{f:tube}
\end{figure}
\vspace{-2em}
\end{center}

Let us consider sphere-packing types $5/4/h16$ and $5/4/h29$. As mentioned by Sowa and Koch (2005) \cite{Sowa2005}, they are built from interconnected `nanotubular' 1-periodic subunits which correspond to square lattice rolled along the vectors (4, 12) and (5, 12), respectively. In formal terms this means that the generating set for $5/4/h16$ and $5/4/h29$ contains a subset referring to a 1-periodic subunit. In the notation of Table~\ref{t:monsters}, in both cases this subset is $\{a, c, d\}$. We have the following presentations:

\begin{center}
\begin{tabular}{l}
$p6_{2}(22) \cong D_{\infty} \times C_2 = \langle a, c, d | a^2, c^2, (ad)^2, (cd)^2, d^{-1}(ac)^2d^{-3}(ac)^4 \rangle$ \\
$p6_{1}(22) \cong D_{\infty} = \langle a, c, d | a^2, c^2, (ad)^2, (cd)^2, adcacd^{-3}(ac)^4d^{-1} \rangle$ \\
\end{tabular}
\end{center}

\noindent Besides relators of length 4 that are obviously inherited from the cycles of the parent square lattice, we immediately notice the relators of length 16 and 17 which correspond to strong rings along the perimeter of a cylinder (Figure~\ref{f:tube}). In connection with embeddings of graphs on surfaces such rings are termed \emph{collar}. The number of such strong rings is equal to the number of geodesics from vertex (0, 0) to vertex (4, 12) respectively (5, 12) in a square lattice which are equal to 1820 and 6188, as confirmed by an independent computation. Any of these \emph{long} geodesics can be taken as a defining relator, and moreover, all corresponding rings are retained intact in both $5/4/h16$ and $5/4/h29$ (\emph{cf.} Schl\"afli symbols in Table~\ref{t:monsters}). Therefore, the extremely high number of strong rings is not a signature for the overwhelming `complexity' of a structure. In a certain sense it can be said that the structure of a Cayley graph appears to be more complicated than a set of rules (relators of a group presentation) by which it is defined. In our opinion, this interpretation holds for all graphs from Table~\ref{t:monsters}.

Let us discuss an exceptional example of another kind, namely when group presentations can be used to find a good upper bound for the size of strong rings in the absense of any \emph{a priori} information. Consider two generating sets $\{a, c, d\}$ and $\{b, c, d\}$ for space group $Pnna$:

\begin{center}
\begin{tabular}{l}
$a=2 (1/4, 0, z): 1/2-x, -y, z$; \\
$b=2 (1/4, 1/2, z): 1/2-x, 1-y, z$; \\
$c=2 (x, 3/4, 1/4): x, 3/2-y, 1/2-z$; \\
$d=\bar1 (0, -1/2, 0): -x, -1-y, -z$.
\end{tabular}
\end{center}

\noindent Geometrically Cayley graphs for both generating sets are built up from helices running along the 2-fold screws at $1/4, y, 1/4$ which are interconnected by edges corresponding to inversion centers. However, in the first case helices are triple ($\langle a, c \rangle$), in the second -- single ($\langle b, c \rangle$). Computed presentations contain relators of length at most 40. Both graphs are isomorphic to each other within the balls of radius~19 (\emph{cf.} \cite{Bab20}), from the 20th sphere onwards coordination sequences become different. A brute-force computation \cite{Zoub24} confirms the presence of such large strong rings:

\begin{center}
\begin{tabular}{l}
$\langle a, c, d \rangle: 10^5.14^{14}.40^{11800}$ \\
$\langle b, c, d \rangle: 10^5.14^{14}.40^{1200}$
\end{tabular}
\end{center}

\begin{center}
\begin{figure}[ht]
\centering
\includegraphics{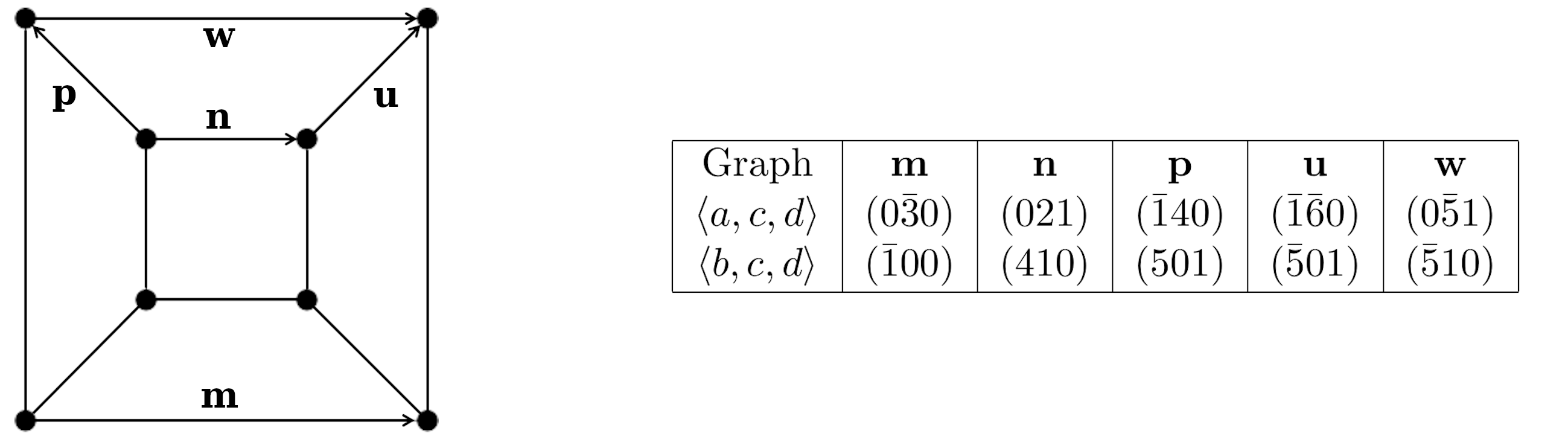}
\caption{Examples of 3-periodic quotient graphs of the 5-periodic cubic graph.}
\label{f:CUBE}
\end{figure}
\end{center}

\noindent A thorough analysis shows that both graphs are quotients of the same parent 5-periodic \emph{minimal} graph by suitably chosen translation groups of rank 2 (Figure~\ref{f:CUBE}). This trivalent 5-periodic graph is vertex-, edge- and strong-ring-transitive, with local configuration of rings characterized by the symbol $14^{21}$. The two translations that are factored out correspond to specific geodesics of length 10 and 40. The number of geodesics of length 10 (they are the same for both quotients) is five that can be easily verified as they correspond one-to-one to the number of strong rings (sharing a vertex) which are formed upon `rolling'. Long `parent' geodesics of length 40 are mostly cut out upon rolling for obvious reasons.

\subsection{4-periodic graphs and more}

Our results for 4-periodic minimal graphs show that ring-transitivity in $\R^4$ is a very rare phenomenon for this class of graphs. Perhaps this fact points towards a general observation that in $\R^4$ highly symmetric structures are not necessarily those of highest translational symmetry (\emph{i.e.}, with a small number of vertex orbits mod $\Z^n$).  While in $\R^3$ the \emph{four} out of 15 minimal graphs are vertex- and ring-transitive ({\bf srs}, {\bf ths}, {\bf dia}, {\bf pcu}), in $\R^4$ only hypercubic lattice and hyperdiamond (out of 111 minimal graphs in total) share this property. Note that these types can be generalized to any dimension \cite{OK91}. 

A few graphs from Tables \ref{t:table1} and \ref{t:table5} are related to each other: 6(3)4 and {\bf hcb} (the latter is the quotient of the former by some rank-2 translation subgroup), {\bf nbo} and {\bf qtz} are quotients of 3(4)1 by some rank-1 translation subgroups \cite{Eon2011}. The translations in 3(4)1 to be factored out correspond to geodesics of length 6 in both cases. However, the number of these geodesics is different (8 vs. 6) that manifests itself in the different number of strong 6-rings per a vertex in {\bf nbo} and {\bf qtz} (Table~\ref{t:table1}). They are indeed the only vertex-transitive 3-periodic quotients of 3(4)1 that can be obtained by factoring out translations corresponding to any geodesic of length 6. Interestingly, geodesics of length 8 do not give rise to any vertex-transitive quotients, those of length $\geq 9$ generate too dense graphs that show self-catenation (see supporting information).

Let us demonstrate interrelations between graphs of different periodicity and their vertex-transitive groups by considering $n$-dimensional diamond (\emph{n}-{\bf dia} in the following), $n\geq2$. Vertices of \emph{n}-{\bf dia} form a \emph{bilattice}\footnote{Two vertex orbits $mod \; \Z^n$ -- hence the term \emph{bilattice}.}, vertex figure corresponds to a regular simplex. The most obvious vertex-transitive group is generated by $(n+1)$ inversions through edge-midpoints. In $\R^2$ and $\R^3$ these are the familiar groups $p2$ ({\bf hcb}) and $P\bar1$ ({\bf dia}), respectively (\emph{cf.} Table~\ref{t:table1}). Let us introduce in $\R^n$ an (affine) coordinate system of unit vectors ($\bf x_1, x_2, ..., x_n$) together with the origin $\bf0$. Inversions through $(n+1)$ points at $\bf 0$, ${\bf{x_1}}/2, {\bf{x_2}}/2, ..., {\bf{x_{n-1}}}/2$ and $1/2 {\bf(x_1 + ... + x_n)}$ generate the desired vertex-transitive group for which \emph{n}-{\bf dia} is the Cayley graph.  In abstract terms, we have a generating set of $(n+1)$ involutions which have to be subject to relators corresponding to the (strong) 6-rings in \emph{n}-{\bf dia}. From Table~\ref{t:table1} we have:

\begin{center}
\begin{tabular}{l}
$p2 = \langle a, b, c |  a^2, b^2, c^2, (abc)^2 \rangle$ \\
$P\bar1 = \langle a, b, c, d |  a^2, b^2, c^2, d^2, (abc)^2, (abd)^2, (acd)^2 \rangle$
\end{tabular}
\end{center}

\noindent In four dimensions:

\begin{center}
\begin{tabular}{l}
$\langle a, b, c, d, f |  a^2, b^2, c^2, d^2, f^2, (abc)^2, (abd)^2, (acd)^2, (abf)^2, (acf)^2, (adf)^2 \rangle$
\end{tabular}
\end{center}

\noindent Proceeding by induction, we can see that in general the group contains $n(n-1)/2$ relators which correspond to the commutators of the generators for the translation group ($ab, ac, ad, af, ...$). More explicitly: $[ab, ac] = bacaabac = bacbac = (abc)^2$.

It is interesting to compare the number of orbits of 6-rings and the number of relators. To count the orbits of 6-rings, it is easier to start with $Aut(n\text{-}{\bf dia})$. Vertex stabilizer in $Aut(n\text{-}{\bf dia})$ is isomorphic to $S_{n+1}$ whereas the full group is generated by appending an inversion through one of the edge midpoints to $S_{n+1}$.  The point group of \emph{n}-{\bf dia} is isomorphic to $S_{n+1} \times C_2$ of order $2(n+1)!$, the stabilizer of a 6-ring is isomorphic to $(D_{3} \times C_2) \times S_{n-2}$. The number of 6-ring-orbits $mod \; \Z^n$ is therefore $n(n-1)(n+1)/6$. The number of 6-rings sharing a common vertex is

\begin{center}
\begin{tabular}{l}
$\dfrac{n(n-1)(n+1)}{6 \cdot 2} \times 6 = \dfrac{n(n-1)(n+1)}{2}$
\end{tabular}
\end{center}

\noindent -- the formula given by O'Keeffe \cite{OK91} without proof. 

\begin{center}
\begin{figure}[ht]
\centering
\includegraphics{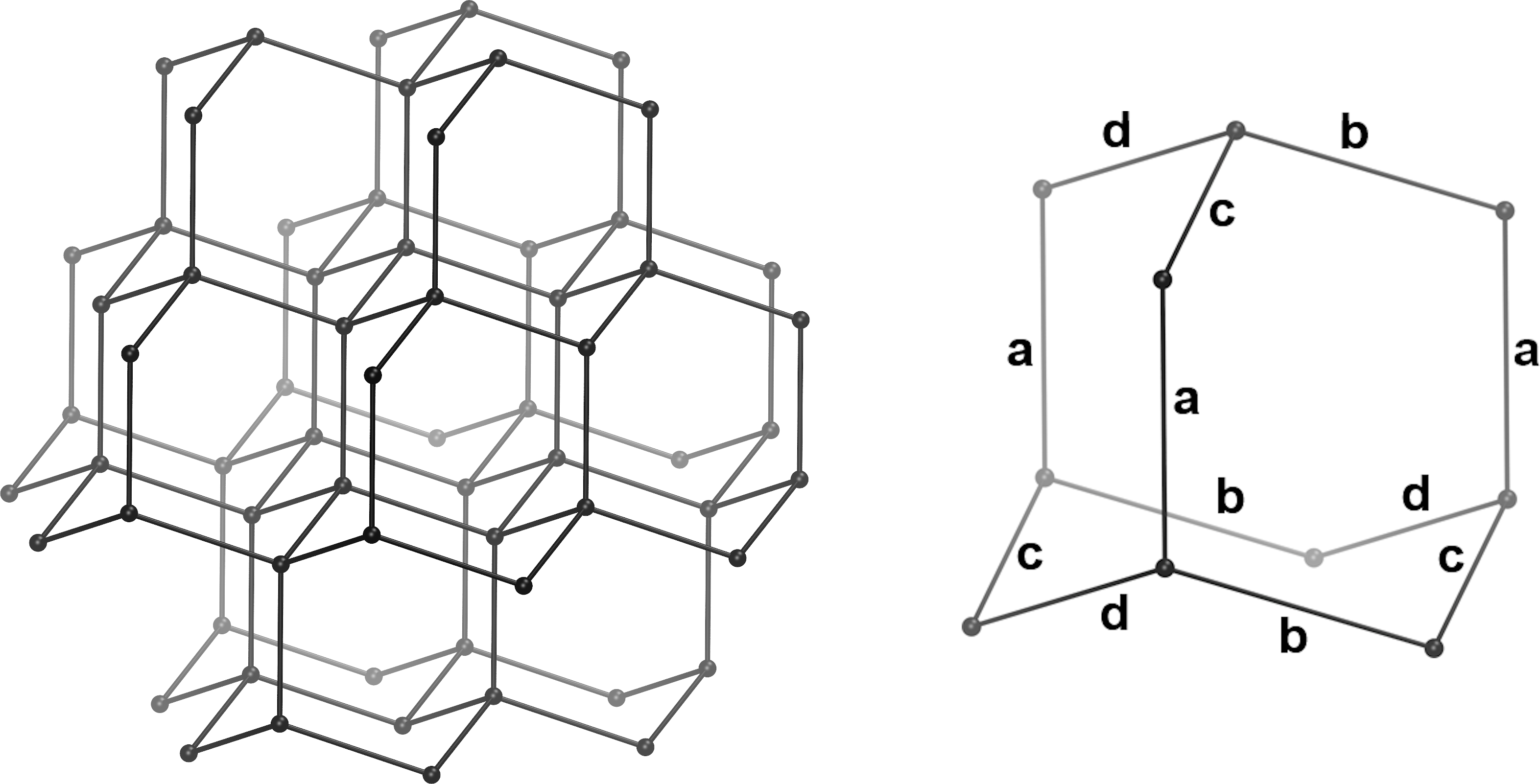}
\caption{Left: a finite portion of the {\bf dia} graph. Right: the adamantane cage as a space-filling polyhedron for space group $P\bar1$. For the edge labels see text.}
\label{f:DIATILE}
\end{figure}
\end{center}

In the group generated by inversions through edge mid-points the orbit of 6-rings in  $Aut(n\text{-}{\bf dia})$ splits up into $n(n-1)(n+1)/6$ orbits (since each 6-ring $mod \; \Z^n$ is stabilized by its own inversion center) whereas the number of relators of length 6 is $n(n-1)/2$.  Only in two dimensions do these numbers coincide. In three dimensions there are 3 relators and 4 orbits of 6-rings. That the fourth orbit does not introduce an independent relator can be seen from the fact that {\bf dia} is an \emph{isohedral} tiling of $\R^3$  by generalized tetrahedra with 6-rings as 2-faces. The cycle space of a 3D tetrahedron has dimension 3 and is generated by any three faces meeting at a vertex. For example, from Figure~\ref{f:DIATILE} it should be clear that the relator $(bcd)^2$ corresponding to the base-face of a tetrahedron can be obtained from another three by a sequence of concatenations and cancellations. However, in higher dimensions this transparent geometrical interpretation is lost.

\section{Conclusion and outlook}

In this work we applied the machinery of computational group theory to study presentations of crystallographic groups and the associated Cayley graphs, especially their `cycle structure'. The obtained presentations can be used to explore subgroup structure of groups, to compute generating functions for coordination sequences using the available tools for automatic groups \cite{Ep92}, incl. the powerful Knuth-Bendix method \emph{etc.} In the supporting information we provide generating functions for coordination sequences of vertex-transitive zeolites which were computed using the KBMag progam \cite{KBMAG}. This computation gives an alternative proof that coordination sequences for certain zeolites coincide (pairs LTA/RHO and ABW/ATN) -- the point that has been puzzling crystallographers for some time \cite{Grosse96}.

As an extension towards more interesting physical applications, we mention here the problem of how the global structure of a vertex-transitive graph is determined by its local structure (that is, by the combinatorial structure of `small' balls). Vertex-transitive graphs occurring in crystal structures \cite{Bab08} are characteristic in a sense that many of them are distinguished by balls of small radii, \emph{i.e.}, (non-isomorphic) graphs which are locally isomorphic within a ball of radius 3 or 4 are rare. This enables finite balls (or closely related Schl\"afli clusters) to be used as isomorphism invariants for infinite graphs, for example, to recognize certain structural motifs in the atomistic models of amorphous materials \cite{Treacy2012}. However, it is not clear whether this \emph{local distinctness} has (mainly) a combinatorial or physico-chemical nature. To partially answer this question, we search for an approach that would allow

\begin{enumerate}
\setlength\itemsep{-0.3em}
\item[(a)] given a Cayley graph a crystallographic group, determine the smallest ball that can be uniquely extended to it;
\item[(b)]given a finite ball, derive all pairwise non-isomorphic Cayley graphs which could be `grown' from it.
\end{enumerate}

Technically the solution could be provided by an algorithm that -- given a finite ball cut out from a Cayley graph of a crystallographic group -- enumerates group presentations corresponding to crystallographic groups by considering all possible labellings of the ball edges by (new) generators. The radius of this ball should be at least $[r/2]$ where $r$ is the length of a longest relator in a short presentation (for which the Cayley graph was originally constructed). For example, for a generating set of four involutions, the ball of radius 3 from the {\bf dia} graph would necessarily produce labellings corresponding to the crystallographic groups from Table~\ref{t:table1}: $P\bar1$, $P2/c$, $C2/c$, $P222_1$, $Pnna$.

In general, however, this problem is technically hard and depends mostly on fast tools for recognizing crystallographic groups of a given dimension from a finite presentation. This recognition problem also arises in the Delaney--Dress tiling theory where it is solved practically by making use of (among other things) the orbifolds of crystallographic groups \cite{Delgado2001} since from the very beginning tilings of simply-connected spaces are considered. This aspect is missing in our purely combinatorial problem that makes it harder to tackle.

\bigskip

\noindent {\bf Acknowledgement}. Part of this work was done in 2022--2023 when the author served as a substitute professor (Vertretungsprofessor) at the Ludwig-Maximilians-Universit\"at M\"unchen, Sektion Kristallographie. Funding within the Hightech Agenda Bayern is gratefully acknowledged.

\bigskip

\noindent {\bf Supporting information}: the $\GAP$ code for computing short presentations, generating sets and presentations for regular subgroups of $Aut({\bf pcu})$, coordinates for graphs from Tables \ref{t:table1} and \ref{t:table2}, generating functions for coordination sequences of zeolites (Tables \ref{t:table2}, \ref{t:table3}), the data on strong rings for 1059 sphere-packing graphs. All the materials can be obtained directly from the author.

\end{document}